\title{ ~~\\ On the distribution of the order over residue classes}
\author{Pieter Moree}
\documentclass[12pt]{article}
\usepackage{amssymb, latexsym, amsfonts}
\def\@ptsize{2}
\newtheorem{Thm}{Theorem}

\newcommand{\qed}{\hfill $\Box$}

\begin{document}
\date{}
\maketitle
{\def\thefootnote{}
\footnote{\noindent 
Max-Planck-Institut f\"ur Mathematik, Vivatsgasse 7,
D-53111 Bonn, Deutschland, E-mail:
moree@mpim-bonn.mpg.de\\}
{\def\thefootnote{}
\footnote{{\it Mathematics Subject Classification (2000)}.
primary 11N37, 11R45, secondary 11N69}}

\begin{abstract}
\noindent For a fixed rational number $g\not\in \{-1,0,1\}$ and
integers $a$ and $d$ we consider the 
set  $N_g(a,d)$ of primes $p$ such that the order of
$g$ modulo $p$ is congruent to $a({\rm mod~}d)$. 
Under the
Generalized Riemann Hypothesis (GRH), it can be
shown that the set $N_g(a,d)$ has a natural
density $\delta_g(a,d)$. Arithmetical properties of $\delta_g(a,d)$ are
described and $\delta_g(a,d)$ is compared with $\delta(a,d)$: the average density
of elements in a field of prime characteristic having order congruent to $a({\rm mod~}d)$.
It transpires that $\delta_g(a,d)$ has a strong tendency to be equal to $\delta(a,d)$, or 
at least to be close to it.
\end{abstract}
\section{Introduction}
\label{introdu}
Let $g\not\in \{-1,0,1\}$ be a rational number. For a rational number $u$, let
$\nu_p(u)$ denote the exponent of $p$ in the canonical factorisation of
$u$ (throughout the letter $p$ will be used to indicate prime numbers).
If $\nu_p(g)=0$, then there exists a smallest positive integer $k$ such
that $g^k\equiv 1({\rm mod~}p)$. We put ord$_p(g)=k$. The number $k$ is
the {\it (residual) order} of $g({\rm mod~}p)$. We let $N_g(a_1,d_1;a_2,d_2)$ 
be the set of primes $p$ with $\nu_p(g)=0$, $p\equiv a_1({\rm mod~}d_1)$ and
ord$_p(g)\equiv a_2({\rm mod~}d_2)$. By $N_g(a_1,d_1;a_2,d_2)(x)$ we
denote the number of primes $p\le x$ in $N_g(a_1,d_1;a_2,d_2)$.
For convenience $N_g(0,1;a,d)(x)$ is denoted as $N_g(a,d)(x)$. By GRH we denote the Generalized Riemann
Hypothesis. By $(a,b)$ and $[a,b]$ we denote the greatest common divisor, respectively 
lowest common multiple of $a$ and $b$.
\begin{Thm} {\rm \cite{Moree2}. (GRH)}. The density $\delta_g(a_1,d_1;a_2,d_2)$ of
the set of primes $N_g(a_1,d_1;a_2,d_2)$ exists. Moreover,
$$N_g(a_1,d_1;a_2,d_2)(x)=\delta_g(a_1,d_1;a_2,d_2){x\over \log x}
+O_{g,d}\left({x\over \log^{3/2}x}\right).$$
\end{Thm}
Our primary interest is in $\delta_g(a,d)=\delta_g(0,1;a,d)$, but in studying this
quantity it turns out to be fruitful to consider $N_g(a_1,d_1;a_2,d_2)(x)$.
Theorem \ref{d=4} for example is obtained from Theorem \ref{4opsplits}.
By $K_{s,r}$ (with $r|s$) we denote the number field $\mathbb Q(\zeta_s,g^{1/r})$, 
where $\zeta_s=\exp(2\pi i/s)$. The density 
$\delta_g(a_1,d_1;a_2,d_2)$ can be expressed in terms of the degrees $[K_{s,r}:\mathbb Q]$ and certain
intersection coefficients (with $(b,f)=1$),
$$c_g(b,f,v)=\cases{1 &if $\sigma_b|_{\mathbb Q(\zeta_f)\cap K_{v,v}}$=~identity;\cr
0 &otherwise,}$$
where $\sigma_b$ is the automorphism of $\mathbb Q(\zeta_f)$ that sends
$\zeta_f$ to $\zeta_f^b$. 
\begin{Thm} {\rm \cite{Moree2}. (GRH)}. 
\label{tweetje}
We have
\begin{equation}
\label{centraal}
\delta_g(a,d)=
\sum_{t=1\atop (1+ta,d)=1}^{\infty}
\sum_{n=1\atop (n,d)|a}^{\infty}{\mu(n)c_g(1+ta,dt,nt)\over 
[K_{[d,n]t,nt}:\mathbb Q]}.
\end{equation}
\end{Thm}
In the proofs of various 
results below the determination of the intersection coefficients $c_g(b,f,v)$ plays an
important r\^ole. Obviously $\mathbb Q(\zeta_f)\cap K_{v,v}$ is a subfield of the maximal
abelian subfield, $K_{v,v}^{\rm ab}$, of $K_{v,v}$. It turns out, see \cite{Moree3}, that
$K_{v,v}^{\rm ab}$ is of the form $\mathbb Q(\zeta_v,\sqrt{\gamma})$ or of the form
$\mathbb Q(\zeta_v,\zeta_{2v}\sqrt{\gamma})$ for some integer $\gamma$ that can
be explicitly given, where the latter
case does not arise if $g>0$. Of course, the action of $\sigma_b$ on $\mathbb Q(\zeta_f)\cap K_{v,v}
=\mathbb Q(\zeta_f)\cap K_{v,v}^{\rm ab}$, which equals $\mathbb Q(\zeta_f)\cap 
\mathbb Q(\zeta_v,\sqrt{\gamma})$
or  $\mathbb Q(\zeta_f)\cap \mathbb Q(\zeta_v,\zeta_{2v}\sqrt{\gamma})$, is easily determined.\\
\indent The distribution of the order over congruence
classes in case $d\nmid a$ seems to have been first studied by
Chinen and Murata \cite{CM} for  $d=4$. In case $d|a$ the problem is much easier and
unconditional results have been obtained, see \cite{M-old, Wiertelak1, Wiertelak2}. In this
case the density is always a rational number. 
Chinen and Murata restricted themselves to the case where $g$ is positive and
not a power of an integer. 
In their method $\delta_g(a,4)$ (for $a$ is odd) is
initially expressed as the sum of two fourfold sums. 
On making everything explicit, they obtained a long formula (distinguishing six cases) for
$\delta_g(a,4)$ which was subsequently simplified by Zagier \cite{Z}.
The author expressed $\delta_g(a,4)$ as a single sum (see Theorem \ref{4opsplits}), which on evaluation
gives a compact formula for $\delta_g(a,4)$ similar to Zagier's. In this
formula an Euler product $A_{\psi}$ appears. We put, for any Dirichlet character $\chi$,
$$A_{\chi}=\prod_{p\atop \chi(p)\ne 0}\left(1+{[\chi(p)-1]p\over [p^2-\chi(p)](p-1)}\right).$$
The constants $A_{\chi}$ turn out to be rather basic in this setting (cf. Theorem \ref{indeks}).
A table of numerical values of $A_{\chi}$, with $\chi$ a Dirichlet character of modulus $\le 12$, is
given in \cite{M-old}.
We let $\cal G$ be the set of rational integers that can not be written as
$-g_0^h$ or $g_0^h$ with $h>1$ an integer and $g_0$ a rational number. Note that
almost all integers $g$ are elements of $\cal G$.
\begin{Thm}
\label{d=4}
{\rm \cite{Moree1}. (GRH).} Let $D(g)$ denote the discriminant of the field $\mathbb Q(\sqrt{g})$.
Let $g\in \cal G$. Then $\delta_g(\pm 1,4)=1/6$ unless
$D(g)$ is divisible by $8$ and has no prime divisor congruent to $1({\rm mod~}4)$, in
which case we have
$$\delta_g(\pm 1,4)=\cases{{1\over 6}\mp {\rm sgn}(g){A_{\psi}\over 8}
\prod_{p|{D(g)\over 8}}{2p\over p^3-p^2-p-1} & if $D(g)\ne \pm 8$;\cr
{7\over 48}\mp {\rm sgn}(g){A_{\psi}\over 8} &if $D(g)=\pm 8$},$$
where $\psi$ denotes the non-trivial character mod $4$.
\end{Thm}
{\tt Remark 1}. We have, on invoking Theorem \ref{twaalf}, $A_{\psi}=0.643650679662525\cdots$.\\
{\tt Remark 2}. Let $\epsilon_1(n)=1$ if $8|n$ and $0$ otherwise. Then we can write,
using e.g., Theorem 2 of Moree \cite{M-old}, for $2\nmid a$:
$$\delta_g(a,4)={\delta_g(1,2)\over 2}+\epsilon_1(D(g)){\rm sgn}(g)A_{\psi}
{(-1)^{a+1\over 2}\over 8}\prod_{p|{D(g)\over 8}}{(1-\psi(p))p\over p^3-p^2-p-1}.$$
An explicit expression for $\delta_g(a,q^s)$ (in terms 
of $A_{\chi}$'s) with $q$ a prime and $g\in \cal G$ is
obtained in \cite{Moree2}, but is omitted here for reasons of space.\\
\indent Theorem \ref{d=4} can be obtained from the following result (with $s=2$):
\begin{Thm} {\rm \cite{Moree1}. (GRH)}. 
\label{4opsplits}
We have, for $a$ odd and $s\ge 1$,
$$N_g(1,2^s;a,4)(x)={\delta_g(1,2^s;1,2)\over 2}{x\over \log x}+O_g({x\over \log^{3/2}x}),$$
and $N_g(3,4;a,4)(x)=\#\{p\le x:p\equiv 3({\rm mod~}4),({g\over p})=1\}/2$
$$+(-1)^{a-1\over 2}{\Delta_g\over 4}{x\over \log x}+O_g({x\over \log^{3/2}x}),$$
where
$$\Delta_g=\sum_{\sqrt{-2}\in K_{2v,2v}\atop 2\nmid v}{h_{\psi}(v)\over [K_{2v,2v}:\mathbb Q]}
-\sum_{\sqrt{2}\in K_{2v,2v}\atop 2\nmid v}{h_{\psi}(v)\over [K_{2v,2v}:\mathbb Q]},$$
with $h_{\psi}$ the Dirichlet convolution of $\psi$ and the M\"obius function, i.e., 
$h_{\psi}(n)=\sum_{d|n}\psi(d)\mu(n/d)$.  
\end{Thm}
Remarkably, despite its arithmetic complexity $\delta_g(3,4;a,4)$ satisfies some
easy properties.
\begin{Thm}
{\rm \cite{Moree1}. (GRH)}.
Write $g=\pm g_0^h$, where $g_0$ is positive and not
an exact power of a rational number.\\ 
{\rm 1)} If $g>0$ and $h$ is even, then $\delta_g(2,3;1,3)\le \delta_g(2,3;2,3)$, otherwise
$\delta_g(2,3;1,3)\ge \delta_g(2,3;2,3)$. We have equality iff
$\mathbb Q(\sqrt{g_0})=\mathbb Q(\sqrt{3})$ and $\nu_2(h)\in \{0,2\}$. The same result
holds with $\delta_g(2,3;*,3)$ replaced by $\delta_g(*,3)$.\\ 
{\rm 2)} If $\delta_g(3,4;3,4)\ne \delta_g(3,4;1,4)$, then
${\rm sgn}(\delta_g(3,4;3,4)-\delta_g(3,4;1,4))={\rm sgn}(g)$.\\
{\rm 3)} If $g\in \cal G$ and $2\nmid a$, then
$\delta_g(3,4;a,4)+\delta_{-g}(3,4;a,4)=1/4.$
\end{Thm}
Let $\delta(p;a,d)$ denote the density of elements in $\mathbb F_p^*$ having order
congruent to $a({\rm mod~}d)$. It is not so difficult to show that the average
density $\delta(a,d)$ of elements of order congruent to $a({\rm mod~}d)$ in a field of
prime characteristic exists. I.e., we have $\lim_{x\rightarrow 
\infty}\sum_{p\le x}\delta(p;a,d)/\pi(x)=\delta(a,d)$, where $\pi(x)$ denotes the number
of primes $p\le x$.
The quantity $\delta(a,d)$ can be studied by fairly elementary methods, but nevertheless turns out to
exhibit behaviour similar to $\delta_g(a,d)$. An interpretation of $\delta(a,d)$ is that it is
the $g$-average of $\delta_g(a,d)$:
\begin{Thm} {\rm \cite{Moree3}. (GRH)}.
We have $${1\over 2x}\sum_{|g|\le x}\delta_g(a,d)=\delta(a,d)+O({1\over \sqrt{x}}).$$
\end{Thm} 
In some sense, if one takes out the Galois
theory and degree aspects of formula (\ref{centraal}), one obtains 
$\delta(a,d)$. More precisely, if one sets $c_g(1+ta,dt,nt)=1$ and
$[K_{[d,n]t,nt}:\mathbb Q]=\varphi([d,n]t)nt$ (this is the maximal degree
possible), then it can be shown that one obtains 
$\delta(a,d)$ out of $\delta_g(a,d)$ (\cite{Moreeaverage}). One has $\delta({\rm odd},4)=1/6$.
It is not difficult to prove that for most integers $g$ with $|g|\le x$ we have that
$D(g)$  has a prime divisor congruent to $1({\rm mod~}4)$. 
(Indeed, the size of the exceptional set is bounded above by $\ll_g x/\sqrt{\log x}$).
Thus from Theorem \ref{d=4} 
we infer that for almost all integers $g$ with $|g|\le x$ we have, on GRH, 
$\delta_g({\rm odd},4)=\delta({\rm odd},4)=1/6$. If $\delta_g(a,4)\ne \delta(a,4)$, then
the difference will be small in absolute value as is also obvious from 
Theorem \ref{d=4}. It turns out that these phenomena hold true in general.
In case $d$ equals a prime power it is still possible to write down an explicit formula
for the density $\delta_g(a,d)$ from which the latter two properties can be similarly inferred 
\cite{Moree2}.
For general $d$ this seems to be difficult. Nevertheless, the following
two results can be proved 
(where $k(d)=\prod_{p|d}p$ is the squarefree kernel of $d$):  
\begin{Thm} {\rm \cite{Moree3}. (GRH)}.
\label{favoriet}
Let $d$ be fixed. There are at most $O_d(x\log^{-1/\varphi(k_1(d))}x)$ integers $g$ with $|g|\le x$ such 
that $\delta_g(a,d)\ne \delta(a,d)$ for
some integer $a$. In particular,
$$(\delta_g(0,d),\dots,\delta_g(d-1,d))=(\delta(0,d),\dots,\delta(d-1,d))$$
for almost all integers $g$, where $$k_1(d)=\cases{k(d) & if $d$ is odd;\cr
4k(d) & otherwise,}{\rm ~and~}k_2(d)=\cases{k(d) & if $d$ is odd;\cr
(4,d/2)k(d) & otherwise.}$$
\end{Thm}
\begin{Thm} {\rm \cite{Moree3}. (GRH)}. 
\label{stellingeen}
Suppose that $g\in \cal G$. Set $D_1=|D(g)/(D(g),d)|$. Then
$$\Big|\delta_g(a,d)-\delta(a,d)\Big|<{3\cdot 2^{\omega(D_1)+2}\over \varphi(D_1)D_1},$$
where $\omega(n)$ denotes the number of distinct prime divisors of $n$.
\end{Thm} 
The following basic result reduces the study of $\delta_g(a,d)$ to that
of $\delta_g(a,k_2(d))$.
\begin{Thm} {\rm \cite{Moree3}. (GRH)}.
\label{peelof}$~$\\
{\rm 1)} If $q$ is an odd prime dividing $d_1$, then $\delta_g(a,qd_1)=\delta_g(a,d_1)/q$.\\
{\rm 2)} If $8|d_1$, then $\delta_g(a,2d_1)=\delta_g(a,d_1)/2$.\\
That is, we have $\delta_g(a,d)=\delta_g(a,k_2(d))k_2(d)/d$.
\end{Thm}
It is easy to see that $\delta(a,d)$ satisfies a similar and slightly stronger
property: $\delta(a,d)=\delta(a,k(d))k(d)/d$.\\
\indent In Theorem \ref{4opsplits} it is seen that there is a difference in behaviour of ord$_p(g)$ when
$p$ is restricted to those primes with $p\equiv 1({\rm mod~}4)$, respectively $p\equiv 3({\rm mod~}4)$.
A similar phenomenon (having a Galois theoretic explanation), is seen to hold in general.
\begin{Thm} {\rm \cite{Moree3}. (GRH)}. Suppose that $(a,d)=(b,d)=1$.\\
{\rm 1)} If $d$ is odd, then
$\delta_g(1,k(d);a,d)=\delta_g(1,k(d);b,d)$.\\
{\rm 2)} If $d$ is even, then
$\delta_g(1,2k(d);a,d)=\delta_g(1,2k(d);b,d)$.
\end{Thm}
On the other hand, if $(a,d)\ne (b,d)$ then it seems that rarely $\delta_g(a,d)=\delta_g(b,d)$, cf. 
Theorem \ref{d=4} with the first part of Theorem \ref{4opsplits}.

\section{On the computation of $\delta_g(a,d)$}
As in Theorem \ref{d=4}, in general $\delta_g(a,d)$ can be expressed in terms of linear
combinations of the constants
$A_{\chi}$ with coefficients coming from certain cyclotomic fields:
\begin{Thm} {\rm \cite{Moree3}. (GRH)}.
\label{indeks}
Let $a$ and $d$ be arbitrary natural numbers. Then there exists an integer $d_1|k_1(d)$ such that
$$\delta_g(a,d)=\sum_{\chi\in G_{d_1}}c_{\chi}A_{\chi}{\rm ~with~}c_{\chi}\in \mathbb Q(\zeta_{o_{\chi}}),$$
where $c_{\chi}$ can be explicitly computed, $G_{d_1}$ denotes the group of Dirichlet characters
modulo $d_1$ and $o_{\chi}$ the order of $\chi$ in $G_{d_1}$.
\end{Thm} 
The following result allows one to evaluate the constants $A_{\chi}$ easily
with ten decimal digit precision and hence, by Theorem \ref{indeks}, the density
$\delta_g(a,d)$.
\begin{Thm} 
\label{twaalf}
{\rm \cite{Moreeaverage}}. Let $p_1(=2), p_2,\ldots$ be the sequence of consecutive primes. Let $\chi$
be any Dirichlet character and $n\ge 31$ (hence $p_n\ge 127$). Then
$$B_{\chi}=AL(2,\chi)L(3,\chi)L(4,\chi)S(n)R_1,$$
$$S(n)=\prod_{k=1}^n\left(1+{\chi(p_k)\over p_k(p_k^2-p_k-1)}\right)
\left(1-{\chi(p_k)\over p_k^3}\right)\left(1-{\chi(p_k)\over p_k^4}\right),$$
$$B_{\chi}=\prod_{p}\left(1+{[\chi(p)-1]p\over [p^2-\chi(p)](p-1)}\right)=
A_{\chi}\prod_{p|d}\left(1-{1\over p(p-1)}\right),$$
$$A=\prod_p\left(1-{1\over p(p-1)}\right)=0.3739558136\cdots,~{\rm and~}{1\over 
1+p_{n+1}^{-3.85}}\le |R_1|\le 1+{1\over p_{n+1}^{3.85}}.$$
\end{Thm}
The factor $AL(2,\chi)L(3,\chi)L(4,\chi)$ in the latter result is the beginning of an
expansion of $B_{\chi}$ in terms of special values of $L$-series:
\begin{Thm} {\rm \cite{Moreeaverage}}.
One has
$$B_{\chi}=A{L(2,\chi)L(3,\chi)\over L(6,\chi^2)}\prod_{r=1}^{\infty}\prod_{k=3r+1}^{\infty}
L(k,\chi^r)^{\lambda(k,r)},~{\rm with~}\lambda(k,r)\in \mathbb Z.$$
\end{Thm}
This formula can be used to approximate $B_{\chi}$ 
(and thus $A_{\chi}$) with even higher numerical precision.
The integers $\lambda(k,r)$ are related to so-called convoluted Fibonacci numbers, 
see \cite{MoreeF}, and exhibit certain monotonicity properties in both the
$k$ and $r$ direction (\cite{MoreeF}). These monotonicity properties are valid in a much
more general setting, see \cite{MoreeW}.

\section{On similar results for the index} 
The index, $[(\Bbb Z/p\Bbb Z)^*:\langle g({\rm mod~}p)\rangle]$, of the subgroup 
generated by $g({\rm mod~}p)$ inside
the multiplicative group of residues mod $p$, is
denoted by $r_p(g)$ and called the {\it (residual) index} mod $p$ of $g$.
For this quantity similar questions can be asked with ord$_p(g)$ replaced by $r_p(g)$.
The results under this replacement turn out to be rather similar, see
\cite{Moree1, Moree2, Moree3, Pappalardi}, however, they are much easier
to establish. A reason for this is that in the latter case
intersection coefficients do not appear. It was in this context that the 
constants $A_{\chi}$ were
introduced by Pappalardi \cite{Pappalardi}. 

\section{On the proofs of the results}
For reasons of space we can only give a small sample here. We sketch the proof of 
Theorem \ref{tweetje}.\\
{\it Sketch of proof of Theorem} \ref{tweetje}. On noting that $r_p(g){\rm ord}_p(g)=p-1$ we obtain that
$N_g(a,d)(x)=\sum_{t=1}^{\infty}V_g(a,d;t)(x)$, where
$$V_g(a,d;t)(x):=\#\{p\le x~:~r_p(g)=t,~p\equiv 1+ta({\rm mod~}dt)\}.$$
In this infinite
sum the terms with $t\ge \sqrt{\log x}$  are less easily individually computed, but since
they are small they can be taken together to form an error term, which
can be estimated by $O(x\log^{-3/2}x)$. If $(1+ta,d)>1$, then there
is at most one prime counted by $V_g(a,d;t)(x)$ and this prime has to divide $d$.
In this way one obtains that 
\begin{equation}
\label{hetgaatmaardoor}
N_g(a,d)(x)=\sum_{t\le \sqrt{\log x},~(1+ta,d)=1}V_g(a,d;t)(x)
+O({x\over \log^{3/2}x}).
\end{equation}
Note that $V_g(a,d;1)(x)$ counts the number of primes $p\equiv 1+a({\rm mod~}d)$ such
that, moreover, $g$ is a primitive root modulo $p$. This function, and indeed  $V_g(a,d;t)(x)$, can
be estimated by a variation of Hooley's classical argument \cite{H}. However, we need to carry this out
with a certain uniformity in $t$ which forces us to keep track of the dependence on $t$ of the
various estimates. Furthermore, as will be explained shortly,  the
additional condition $p\equiv 1+ta({\rm mod~}dt)$ is responsible for bringing in the Galois
theoretic intersection coefficients $c_g(1+ta,dt,nt)$. By inclusion and exclusion we
find that 
\begin{equation}
\label{pufpuf}
V_g(a,d;t)(x)=\sum_{n=1}^{\infty}\mu(n)\#\{p\le x:p\equiv 1+ta({\rm mod~}dt),~nt|r_g(p)\}.
\end{equation}
The counting functions in the latter sum can be estimated by an effective form of Chebotarev's
density theorem, cf. Theorem 3
of \cite{Moree1} and the discussion immediately
following that theorem. Namely, we are interested in those primes $p\equiv 1+at({\rm mod~}dt)$ that
split completely in $K_{nt,nt}:=\mathbb Q(\zeta_{nt},g^{1/nt})$. These primes must have
a Frobenius $\sigma$ in $K_{[n,d]t,nt}$ with the property that $\sigma|_{\mathbb Q(\zeta_{dt})}=\sigma_{1+ta}$
and $\sigma|_{K_{nt,nt}}=$id. If such a $\sigma$ exists then certainly we must
have $\sigma_{1+ta}|_{\mathbb Q(\zeta_{dt})\cap K_{nt,nt}}=$id, i.e. $c_g(1+ta,dt,nt)=1$.
Indeed, such a $\sigma$ turns out to exist iff $c_g(1+ta,dt,nt)=1$.
On applying Chebotarev's density theorem one then finds, assuming the Riemann Hypothesis (RH)
holds for the field $K_{[d,n]t,nt}$, that
$$
\#\{p\le x:p\equiv 1+ta({\rm mod~}dt),~nt|r_g(p)\}={c_g(1+ta,dt,nt)\over [K_{[d,n]t,nt}]}
{\rm Li}(x)+O(\sqrt{x}\log x).
$$
Again there is a problem with the tail in the series in (\ref{pufpuf}), but again it can
be reasonably estimated and one obtains that $$V_g(a,d;t)(x)=\sum_{P(n)\le 
(\log x)/6}\mu(n)\#\{p\le x:p\equiv 1+ta({\rm mod~}dt),~nt|r_g(p)\}+E(x),$$
where $P(n)$ denotes the greatest prime factor of $n$ and $E(x)$ is
the estimate for the tail. On combining the latter two
displayed estimates,
one then arrives at a usable estimate for $V_g(a,d;t)(x)$. On combining
this with (\ref{hetgaatmaardoor}), the proof of Theorem \ref{tweetje} is then easily
completed. \qed\\

Working with the sharpest known unconditional version of Chebotarev's density theorem leads to an error
term which is too weak for our purposes. Actually, as is clear from the
above sketch it is not required to assume GRH. It suffices to assume RH for the number fields 
involved in the proof. So for Theorem \ref{tweetje} it suffices 
to require RH for the number fields $K_{[d,n]t,nt}$ with $n$ squarefree, $(n,d)|a$
and $(1+ta,d)=1$.

\section{Numerical experiments}
The problem considered here allows for numerical experiments. We give here a small sample
of data so obtained. In the cases studied, the numerics seemed to agree well with the 
theoretical predictions.\\
$~$\\
\centerline{{\bf Table 1:} Experimental and theoretical densities for $d=5$}
\begin{center}
\begin{tabular}{|c|c|c|c|c|c|}\hline
 $\delta\backslash a$ & $0$ & $1$ & $2$ & $3$ & $4$\\ \hline
$\delta(*,5)$ & $0.20833{\underline 3}$ & $0.23542{\underline 1}$ & $0.17799{\underline 3}$ & $0.23400{\overline 3}$ & $0.14424{\overline 8}$ \\ \hline
$\approx \delta_{-11}(*,5)$ & $0.208347$ & $0.235422$ & $0.178007$ & $0.233974$ & $0.144250$ \\ \hline
$\delta_{-11}(*,5)$ & $\delta(0,5)$ & $\delta(1,5)$ & $\delta(2,5)$ & $\delta(3,5)$ & $\delta(4,5)$ \\ \hline
$\approx \delta_{-5}(*,5)$ & $0.208348$ & $0.264146$ & $0.194858$ & $0.233282$ & $0.099365$ \\ \hline
$\delta_{-5}(*,5)$ & $\delta(0,5)$ & $0.26413{\underline 5}$ & $0.19486{\overline 5}$ & $0.23329{\overline 4}$ & 
$0.09937{\underline 1}$ \\ \hline
$\approx \delta_2(*,5)$ & $0.208333$ & $0.240673$ & $0.178706$ & $0.229270$ & $0.143017$ \\ \hline 
$\delta_{2}(*,5)$ & $\delta(0,5)$ & $0.24068{\overline 1}$ & $0.17869{\overline 1}$ & $0.22926{\underline 4}$ & 
$0.14302{\underline 9}$ \\ \hline
$\approx \delta_5(*,5)$ & $0.208348$ & $0.232581$ & $0.292840$ & $0.054488$ & $0.211742$ \\ \hline
$\delta_5(*,5)$ & $\delta(0,5)$ & $0.23258{\overline 5}$ & $0.29284{\overline 8}$ & $0.05449{\overline 3}$ & $0.21173{\overline 7}$ \\
\hline
\end{tabular}
\end{center}
If an entry is in a row labelled $\approx \delta_g(*,5)$ and in column $a$, then the number given 
equals $N_g(a,5)(x)/\pi(x)$ rounded to 6 decimals with $x=2038074743$ (and hence $\pi(x)=10^8$). The
theoretical values are given with 6 digit precision, with a bar over the last digit indicating that
if the number is to be rounded off, it should be rounded upwards. The density $\delta(0,5)=5/24$ 
(unconditional result).\\

\noindent {\bf Acknowledgements}. Most of my papers mentioned in the references were written
whilst I was working in the PIONEER-group of Prof.~E.~Opdam at the University of
Amsterdam (2000-2004). Several of these papers were completed whilst I was enjoying the
inspiring atmosphere of the Max-Planck-Institute in Bonn. I thank both institutes for
their hospitality. The data given in Table 1 were calculated using a $C^{++}$ program
kindly written by Dr.~Yves Gallot. My special thanks go to Dr.~Paul Tegelaar for his
unfailing support and interest over the years.

\end{document}